# Impact of Service Sector Loads on Renewable Resource Integration

Nina Voulis[a,*], Martijn Warnier[a], Frances M.T. Brazier[a]

[a]*Section Systems Engineering, Faculty of Technology, Policy and Management, Delft University of Technology, Jaffalaan 5, 2628 BX Delft, the Netherlands*

**Highlights**

- Improved urban demand modelling using both household and service sector load profiles
- Detailed service sector load profiles based on an extensive data aggregation approach
- Estimation of service sector demand impact in various solar and wind supply scenarios
- New approach for the influence of time and weather on renewable resource potential
- Estimation of service sector demand impact in different time and weather conditions

**Abstract**

Urban areas consist of a mix of households and services, such as offices, shops, schools, etc. Yet most urban energy models only consider household load profiles, omitting the service sector. Realistic assessment of the potential for renewable resource integration in cities requires models that include detailed demand and generation profiles. Detailed generation profiles are available for many resources. Detailed demand profiles, however, are currently only available for households and not for the service sector. This paper addresses this gap. The paper (1) proposes a novel approach to devise synthetic service sector demand profiles based on a combination of a large number of different data sources, and (2) uses these profiles to study the impact of the service sector on the potential for renewable resource integration in urban energy systems, using the Netherlands as a case study. The importance of the service sector is addressed in a broad range of solar and wind generation scenarios, and in specific time and weather conditions (in a single scenario). Results show that including the service sector leads to statistically significantly better estimations of the potential of renewable resource integration in urban areas. In specific time and weather conditions, including the service sector results in estimations that are up to 33% higher than if only households are considered. The results can be used by researchers to improve urban energy systems models, and by decision-makers and practitioners for grid planning, operation and management.

## 1. Introduction

Distributed renewable energy resources are becoming of increasing importance to urban energy systems, posing many new technical, organisational and infrastructural challenges. Addressing these challenges requires a deep understanding of system behaviour at neighbourhood and municipality scales [1, 2, 3, 4]. Urban energy system models can serve this purpose. However, most existing models are based solely on households and do not include service sector load profiles. As real urban areas consist of both households and services[1], such as shops, schools, etc, each with their own energy load profiles, omitting the service sector in energy models is not realistic. The annual demand of the service sector is on par with that of the residential sector in developed countries [7, 8, 9]. However, their load profiles differ considerably [10, 11]. Therefore, urban energy models need to be extended to include the service sector, improving the understanding of the potential of renewable resource integration in cities.

*Corresponding author
  *Email address:* n.voulis@tudelft.nl (Nina Voulis)

This paper proposes a systematic method to devise synthetic load profiles based on a large number of different data sources. Applying the proposed method for the Netherlands, this paper quantifies the impact of the service sector in future urban energy systems with a high penetration of renewables. Three renewable resource integration metrics are compared for a realistic mix of residential and service sector loads, and for residential loads only: (1) mismatch between renewable generation and demand, (2) renewable resource utilisation, and (3) self-consumption. These metrics are first studied for a broad range of solar and wind generation mix scenarios. Second, metrics are compared for different times of the day, days of the week, and weather conditions for a single scenario.

This is the first fundamental study that systematically addresses the impact of service sector loads on renewable resource integration in urban areas. The results of this paper are primarily of interest to researchers in urban energy systems, and to decision-makers and practitioners for grid planning, management and operation, for example, to inform decision-making on storage location, demand response programs and grid reinforcement.

*1.1. Service Sector Demand*

The service sector, also termed the commercial, business or tertiary sector, is comprised of a highly heterogeneous group of energy consumers. Although the many definitions of the sector differ, most include non-manufacturing commercial activities and exclude agriculture and transportation [11, 12]. This paper defines the service sector as the collection of non-manufacturing commercial and governmental activities, excluding agriculture, transportation, power sector, street lighting and waterworks.

The service sector power demand in developed countries currently accounts for one quarter to one third of the total national power demand, and is thus on par with residential demand [7, 8, 9]. Current estimations [8, 13, 14] indicate that in 2050, the demand shares of the service and the residential sectors are projected to increase to 40% each, at the expense of the industry demand, which will account for only 20% of the total national demand.

Despite the importance of the service sector in urban demand, most studies on (future) urban energy system models are based on residential load profiles only. Mikkola and Lund [4] are a notable exception. The authors focus on spatiotemporal modelling of urban areas for energy transition purposes. They include service sector loads in their first case study of Helsinki. The service sector demand profiles for Helsinki are based on German profiles, on the assumption that these profiles are comparable in Finland and Germany. No service sector demand data are referenced by the authors for their second case study of Shanghai. Although a number of other studies consider urban energy systems at the neighbourhood or municipality level, they do not include service sector demand profiles [1, 2, 3]. For instance, Fichera *et al.* [1] study how the integration of distributed renewables in urban areas can be improved using a complex networks approach. Despite considering an entire urban area, the authors model only residential loads in their numerical case study. Hachem [2] describes a neighbourhood designed to increase energy performance and decrease green-house gas emissions and shows that the type of neighbourhood (mixed-use versus residential) has an effect on local renewable energy utilisation. The author models service sector buildings based on a combination of residential data and commercial building code specifications, but does not use detailed service sector profiles. Alhamwi *et al.* [3] focus on the geographical component of urban energy system modelling, and include service sector buildings in their work, yet the authors leave the acquisition of detailed temporal service sector profiles unaddressed.

---

[1]This paper considers only residential and service sectors, not the industrial. Industry is typically located on city edges and has case-specific requirements to transition to renewable generation [5, 6].

These studies illustrate that the importance of the service sector in urban energy systems is increasingly acknowledged, yet that it remains difficulty to obtain detailed service sector demand profiles. The publicly available data on service sector demand are primarily concerned with demand characterisation per subsector (*e.g.*, offices) [11, 15], or per subsector and end-use (*e.g.*, office lighting) [10, 12, 16]. The lack of detailed data is a serious limitation for the assessment of the potential impact of the service sector on renewable resource integration in urban energy systems. This issue is addressed in this paper.

*1.2. Importance of Detailed Demand Data*

Detailed demand data, in combination with detailed generation data, are necessary to realistically assess the impact of interventions designed to increase the integration of renewables in future urban power systems. In traditional power systems, power balance is maintained through generation dispatch, which follows a variable, immutable load [17]. The sizing and operation of dispatchable generation is based on load characteristics such as peak load [18, 19]. Future power systems with a high share of non-dispatchable renewables require different balancing and management approaches, such as demand response and storage. The choice of the best approach or their combination depends on the timing and the extent of mismatches between load and generation. These mismatches depend on (1) the type of load and (renewable) generation, (2) the time (*e.g.*, time of the day, day of the week), and (3) the weather. Understanding the interactions between these factors requires detailed demand and generation data.

On the generation side, detailed profiles are publicly available, or can be constructed from publicly available weather data (*e.g.*, [20] in general, and [21] for the Netherlands). On the demand side, primarily residential profiles are available. These profiles cover only a part of mixed urban demand. Some non-residential profiles are published, including standard load profiles for specific connection types (*e.g.*, [22] for the Netherlands), and country-level load profiles (*e.g.*, [23]). These profiles are not suitable to model mixed urban demand, nor to assess the role of the service sector. The standard load profiles lack essential metadata for non-residential loads, which thwarts their use for the estimation of demand in real urban areas. The country-level load profiles do not give sufficient insights at a local level and can therefore not be used at municipality and neighbourhood scales.

This paper proposes a method to overcome the current demand data scarcity by devising synthetic service sector load profiles, and combining them with the available residential demand profiles to estimate the demand profiles of mixed urban areas. The synthetic load profiles are constructed for an area of interest (the Netherlands, in this paper) and are based on a combination of (1) reference building models of the Unites States (U.S.) Department of Energy [24], and (2) a large number of U.S. and local (Dutch) service sector building use data sources, which are used to scale U.S. reference buildings to the local (Dutch) context.

*1.3. Contributions*

The main contribution of this paper is the systematic assessment of the impact of service sector loads on renewable resource integration in urban areas. This contribution includes:

1. A systematic method to devise detailed synthetic service sector demand profiles based on reference building models and building use data.

2. Quantitative results showing the impact of service sector loads on renewable resource integration metrics for a broad range of renewable resource penetration scenarios.

3. A novel time and weather dependency classification system that enables systematic assessment of metrics that depend on both time and weather.

4. Quantitative results showing the impact of service sector loads on renewable resource integration metrics for different times of the day, days of the week, and weather conditions.

The values for the metrics reported in this paper are specific for the Netherlands. However, the methodology described can be used as a template to assess the impact of the service sector on renewable resource integration in other countries, municipalities, and neighbourhoods. Qualitative conclusions on the impact of the service sector on renewable resource integration in urban areas are assumed to hold for other developed countries as the shape of the service sector load profiles is comparable across countries [10, 24, 25, 26]. The methodology and results are of potential interest to both researchers and practitioners. Improved service sector load profiles can further the research field, for instance through combination with recently developed spatiotemporal urban energy models [1, 3, 4]. Practitioners, such as urban planners, distribution system operators, and aggregators can apply the results for improved grid planning, operation and management.

The remainder of this paper is structured as follows. Section 2 presents the theoretical rationale for explicit consideration of service sector loads. Section 3 describes the systematic approach used to devise synthetic service sector load profiles. Section 4 outlines the methods used for data collection and profile calculation for the Dutch case study. Section 5 provides more details on two renewable resource integration experiments. Section 6 presents the results of these experiments, which are further discussed in Section 7. A final conclusion is given in Section 8.

## 2. Rationale

Urban areas are typically a mix of residential and service sector loads. Residential load profiles are more readily available, making them an attractive proxy for urban areas as a whole. However, residential and service sector load profiles differ considerably. Fig. 1 illustrates the difference in load profiles between (1) residential loads-only and (2) mixed residential and service sector loads. This paper hypothesises that the assessment of renewable resource integration in urban areas based on household demand only is misleading. In particular, it hypothesises that substituting realistic, mixed residential and service sector load profiles with only households load profiles leads to significant misestimations of renewable resource integration metrics. In this section, a theoretical intuition supporting this hypothesis is developed for two metrics: mismatch (between renewable generation and demand), and renewable energy utilisation. These form the basis for the metrics used in the remainder of this paper.

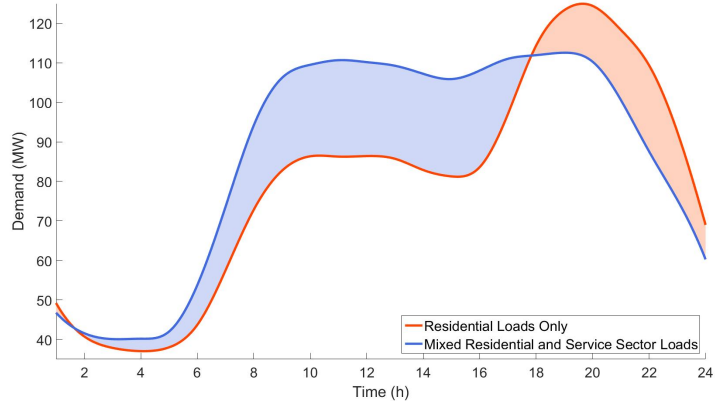

Figure 1: Comparison of load profiles of residential loads-only (orange) and mixed residential and service sector loads (blue) on an average weekday. The shaded area represents the cumulative load difference. The blue shaded area shows energy consumption underestimation by residential loads-only as compared to mixed urban loads. The orange shaded area represents the contrary case. Note that for the mixed load the proportions of residential and service sector are chosen to be representative for the Netherlands.

2.1. Load Type Comparison

Two load types are compared, (1) residential loads-only, denoted by $L_r$, and (2) mixed residential and service sector urban loads, denoted by $L_m$. Let $h(t)$ and $s(t)$ represent respectively household and service sector load over time. Then,

$$L_r(t) = \phi \cdot h(t) \tag{1}$$
$$L_m(t) = h(t) + s(t) \tag{2}$$

Note that $L_r$ is scaled by a factor $\phi$ to ensure that $\sum_{t=1}^{8760} L_r = \sum_{t=1}^{8760} L_m$ for hourly steps of $t$ (with 8760 hours in a year).

2.2. Mismatch

Generation and load must be in perfect balance for the proper operation of the power system. Let mismatch $M$ be the difference between generation $G$ and coinciding load $L$, thus $M = G - L$. For a given generation $G$, the *difference* in mismatch between calculations considering residential and mixed load equals:

$$\begin{aligned}\Delta M &= M_r - M_m \\ &= s(t) - (\phi - 1) \cdot h(t)\end{aligned} \tag{3}$$

From Eq. 3 follows that the difference in mismatch $\Delta M$ only depends on the loads $h(t)$ and $s(t)$, $\Delta M$ does not depend on the generation $G$. Fig. 1 shows that $\Delta M > 0$ during the day and $\Delta M < 0$ in the evening. This suggests that in urban areas with a mixed demand, power imbalance calculations based on residential load only will lead to underestimations of the mismatch between supply and demand during the day and overestimations of this mismatch during the evening. Numerical values and statistical significance of these errors are shown in the following experimental sections.

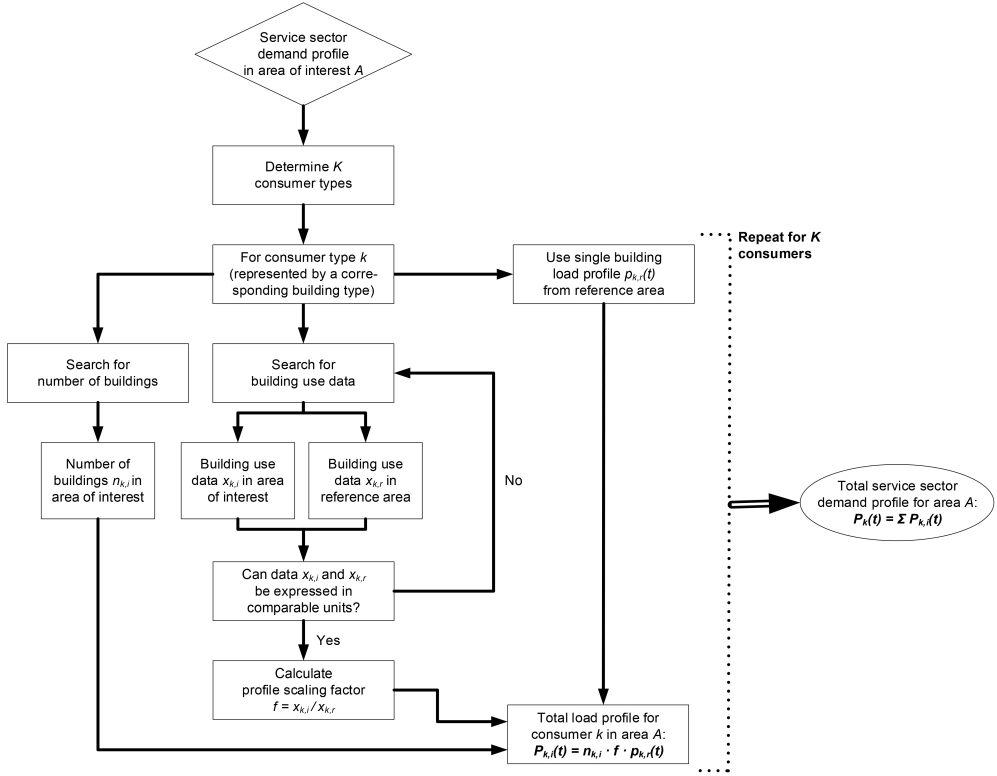

Figure 2: Synthetic service sector profile calculation methodology for an area of interest. Based on load profiles available in a reference area, and building use data from both areas, a detailed service sector load profile for the area of interest can be constructed.

*2.3. Renewable Energy Utilisation*

A similar analysis can be carried out for renewable energy utilisation, denoted by $R$. Assuming no storage or load flexibility, $R$ is computed as:

$$R = \begin{cases} G & \text{if } G \leq L \\ L & \text{if } G > L \end{cases} \quad (4)$$

Given a generation $G$, the *difference* in renewable energy utilisation between residential loads-only and mixed loads equals:

$$\Delta R = R_r - R_m \quad (5)$$

Expanding Eq. 5 yields $\Delta R$ as a function of both renewable generation and load. As both are time and weather dependent, the assessment of renewable energy utilisation requires an analysis of time and weather interactions, in addition to correct load profile estimation. This paper introduces a model which can quantify the influence of load type, both on an annual basis, and in specific time and weather conditions. This model is described in Section 4. First, the general methodology proposed to devise synthetic service sector demand profiles is outlined.

## 3. Methodology: Synthetic Service Sector Load Profiles

Measured detailed service sector load profiles are scarce and rarely available for a specific area of interest. Therefore, this paper proposes a systematic approach to devise detailed synthetic service sector load profiles for an area of interest based on (1) load profiles available in a reference area, and (2) building use data (*e.g.*, occupancy, floor area, etc) from both the area of interest and the reference area. Building use data are used to calculate scaling factors.

The approach depicted in Figure 2 can be summarised as follows. For an area of interest $A$, determine $K$ service sector consumer types. For each consumer type $k$, represented by a building or a collection of buildings (*e.g.*, hospitals), find building use data for both the area of interest ($x_{k,i}$), and the reference area ($x_{k,r}$). Verify that building use data can be expressed in comparable units (*e.g.*, number of beds in a hospital for both areas). Use these data to calculate the profile scaling factor $f = x_{k,i}/x_{k,r}$. The total load profile for the consumer type $k$ from the area of interest equals the reference load profile multiplied by the scaling factor and the number of buildings in the area of interest. Repeat for all $K$ consumer types. Obtain the total service sector demand profile for area $A$ by summing up the profiles of individual consumer types.

This approach can be validated either for each customer type $k$ separately, if annual demand data for each customer type in the area of interest are available, or lumped for all service sector consumers if only total service sector demand data are available.

## 4. Dutch Case Study: Data and Simulation

This paper focuses on the assessment of the impact of service sector loads on renewable resource integration. The Netherlands is chosen as area of interest. The methodology described above is used to devise local synthetic service sector load profiles. These profiles are combined with household load profiles, and solar and wind power generation profiles to create a realistic urban energy model. The influence of load type, and of time and weather on renewable resource integration metrics is studied using a novel simulation model (developed in Matlab [27]). Load type effects are assessed by comparing two load cases: *residential* load only, and *mixed* residential and service sector load. Time and weather effects are studied using a novel time interval classification system. The approach is conceptually shown in Fig. 3. It consists of three steps: (1) data collection, (2) profile modelling, and (3) renewable resource integration experiments. The first two steps, the core of the simulation model, are outlined in this section. The experiments are described in Section 5.

Synthetic load and generation profiles are calculated based on a large number of data sources. To ensure spatial and temporal consistency, all calculations are done for the same area (the Netherlands) and the same period (2014), taking into account official Dutch holidays and daylight saving times. The network is assumed to be a "copper plate". All resulting profiles have an hourly granularity.

### 4.1. Load Types

Two types of loads are defined: *residential load* (comprised of household loads only) and *mixed load* (a mix of household and service sector loads). For both load types, household load is represented by a single average Dutch household profile. For the mixed load type, the service sector load is calculated as a weighted sum of thirteen reference building load profiles (described below in more detail). To ensure that the two load types are comparable, an equal annual cumulative consumption (710 GWh/year) is used for both load types. To achieve this, the residential load is weighted by a factor $\phi$. (In this model, $\phi = 2.03005$, the ratio between the total mixed consumption of 7.10481 GWh and the household consumption for 100 000 households, 3.49982 GWh, see Eq. 1).

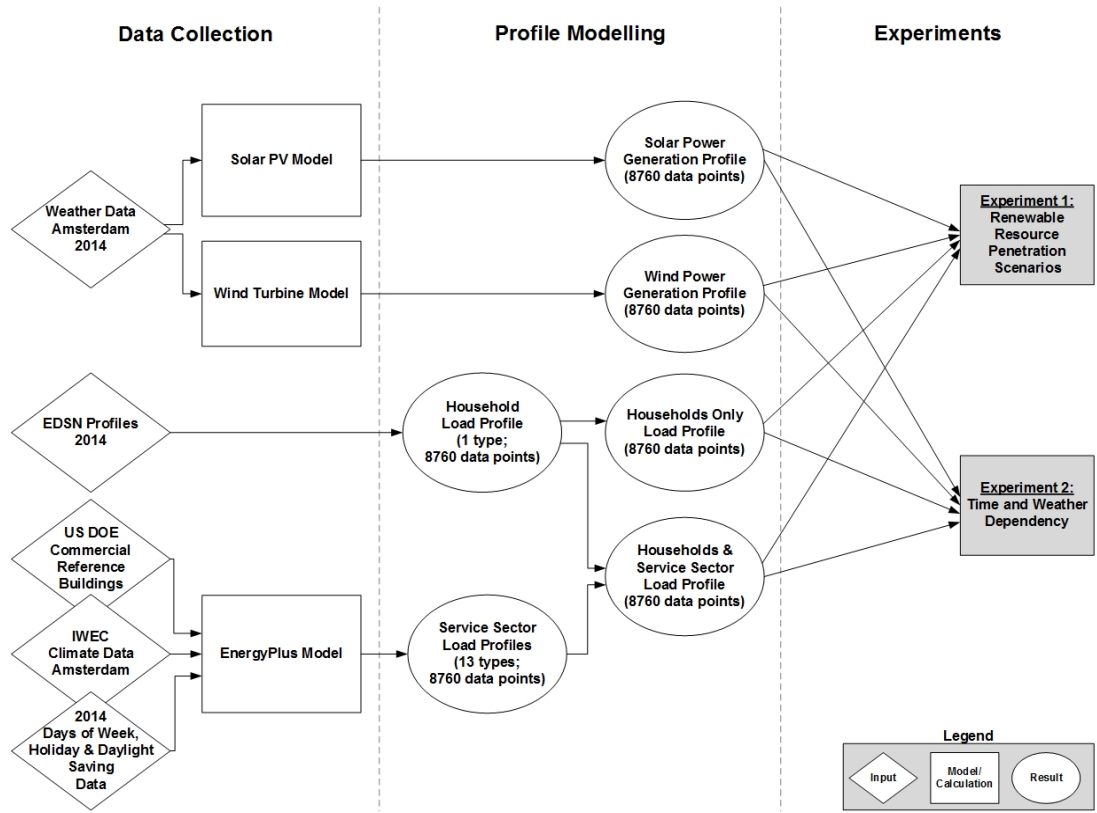

Figure 3: Load and generation profiles calculation. This flow diagram shows how different data sources are combined. The resulting profiles are used in two experiments.

#### 4.1.1. Household Load Profile

Household demand data are obtained from [22]. The average yearly household consumption is assumed to be 3500 kWh [28]. The selected profile describes the *average* Dutch residential load. The use of this single average profile is assumed to be representative at the scale used in the simulations in this paper (100 000 or 203 005 households, depending on load type) since the combined profile of such a large number of similar consumers is expected to regress to the mean profile [29].

#### 4.1.2. Service Sector Load Profiles

Detailed service sector load profiles are modelled based on United States Department of Energy (U.S. DOE) reference building data [24] and a large number of U.S. and Dutch building use data which are used to scale U.S. buildings to the Dutch context (see Section 3 for Methodology, and below and in the Appendix for numerical values). The resulting profiles are combined into a single Dutch service sector profile, which together with the household load profile, is used to model mixed urban loads. The calculations used to scale U.S. reference buildings, as well as the data sources are described in detail in the Appendix (Section A1). An overview is provided in Table 1. The approach can be summarised as follows.

The U.S. DOE publishes reference building data for 16 types of commercial buildings, of which 13 types, relevant for the Netherlands, are used in this paper (see first column in Table 1). The scaling factor for each reference building is based on a comparison of key building use data between U.S. and Dutch buildings (termed "Building Scaling Factor Data" in Table 1). For instance, for hospitals, the average number of patient beds per hospital is used (161 in the U.S. [30] and 316 in the Netherlands [31]). Thus, to model Dutch hospitals using the U.S. reference building and retaining the same service level, the number of Dutch hospitals is multiplied by 2 (316/161 = 2). Further, this study assumes an urban environment of 100 000 households. Thus, the average number of reference buildings (rounded to the nearest integer) per 100 000 households is used (see second column in Table 1). Similar calculations for each reference building type are presented in Section A1.

Table 1: Service sector consumers considered in this paper. A single representative service sector profile is constructed as a weighted sum of the thirteen types of service sector consumers. Column 1: overview of the thirteen service sector consumer types considered. Column 2: estimated number of U.S. reference buildings in the Dutch context per 100 000 households. Column 3: building utilisation data used to calculated the building scaling factors as presented in column 2. Column 4: Roof area available per building for solar panels. Column 5: references used to scale U.S. reference buildings to the Dutch context (for detailed calculation description, see Section A1).

| Service Sector Consumer Type (Reference Buildings) | Number of Reference Buildings (*per 100 000 Households*) | Building Scaling Factor Data | Roof Area Available for Solar PV panels ($m^2$) | References |
|---|---|---|---|---|
| Hospital | 3 | Patient Beds | 4484 | [24, 30, 31] |
| Large Hotel | 1 | Rooms | 1891 | [24, 32, 33] |
| Small Hotel | 16 | Rooms | 1003 | [24, 32, 33] |
| Large Office | 9 | Floor area | 3860 | [24, 34, 35] |
| Medium Office | 47 | Floor area | 1661 | [24, 34, 35] |
| Small Office | 6 | Floor area | 511 | [24, 34, 35] |
| Primary School | 32 | Students | 6871 | [24, 36, 37] |
| Secondary School | 9 | Students | 9796 | [24, 36, 38] |
| Stand Alone Retail | 177 | Floor area | 2294 | [24, 39] |
| Supermarket | 12 | Floor area | 4181 | [24, 40, 41, 42] |
| Restaurant | 170 | Restaurants | 511 | [24, 43] |
| Quick Service Restaurant | 189 | Restaurants | 232 | [24, 43] |
| Warehouse | 163 | Employees | 4835 | [24, 11, 44] |

The service sector profiles themselves are obtained using the DOE EnergyPlus modelling software [45]. This software builds demand profiles based on the building age, climate data, and the building location. As the simulations assume a future situation, new construction (post-2004) standard is used. To create profiles representative for the Netherlands, Amsterdam climate data are used [20]. Finally, the location match in terms of climate zone is based on both the ASHRAE climate classification [46] and the available U.S. locations for the reference models, yielding Seattle as the closest match for Amsterdam. This location match ensures that adequate heating and cooling requirements are taken into account.

*4.2. Generation Profiles*

Solar and wind power generation are modelled using weather data from the Royal Netherlands Meteorological Institute (KNMI) [21].

*4.2.1. Solar Power Generation*

Solar power generation is modelled using the Matlab model developed by Walker [47]. The technical specifications are based on Solarex MSX-60 photovoltaic (PV) panels [48]. This paper assumes that solar panels can be placed on roofs of residential and service sector buildings. The roof area constrains the number of solar panels which can be used. For the service sector, the maximal available roof area is calculated as the ratio between the total floor area and the number of storeys [24] (see fourth column in Table 1). For households, an average roof area of 33 m$^2$ is used [49]. All roofs are assumed to allow for optimal positioning of solar panels. The overestimation of the solar power thus generated, is offset by an underestimation of the solar panel efficiencies, which are rising by 1.0 to 1.2% per year [50]. These two assumptions are expected to balance out between 2030 and 2050, both recurring horizons in literature for scenarios assuming high renewables penetration scenarios (*e.g.* [51, 52]).

*4.2.2. Wind Power Generation*

Wind power generation is modelled as described in [53], for community-size wind turbines of the type 500 kW EWT DIRECTWIND 52/54-500kW [54]. These turbines are 50 m high, have a 54 m rotor diameter and a nominal capacity of 500 kW. The cut-in and cut-out windspeeds of respectively 2.5 m/s and 25 m/s are included in the model. Wind power output is calculated using the following equation [53]:

$$P_{wind} = \frac{1}{2} * \rho * A * V^3 * C_p \qquad (6)$$

where:

$P_{wind}$ : wind power output

$\rho$ : air density, calculated using corresponding temperature and pressure data [21]

$A$ : rotor area, 2290 m$^2$ for modelled turbine

$V$ : wind speed [21], corrected for the height of the wind turbine: $V = V_0 * (H_{turbine}/10)^{0.15}$

$C_p$ : power coefficient, 0.35 for the modelled turbine

*4.3. Service Sector Modelling Validation*

This paper relies on the combination of a large number of openly available data sources to model detailed service sector demand profiles. The best validation for this approach is arguably the comparison of the resulting synthetic profiles with statistically representative, real, measured profiles. However, such profiles are currently not publicly available. That is the very issue this paper is seeking to overcome by estimating service sector profiles and showing the importance of the sector for renewable resource integration. The validation used in this paper thus relies on a different approach, with both a quantitative and a qualitative component.

*4.3.1. Quantitative Validation*

The obtained results are compared with *cumulative* annual Dutch service sector load data, which are openly available but do not suffice to assess the impact of the service sector on renewable resource integration. The Netherlands Environmental Assessment Agency (Planbureau voor de Leefomgeving, PBL) attributes 43.8 TWh of the Dutch annual electricity consumption to the service sector, waste and wastewater treatment, and agriculture and fisheries combined [55]. Solely the service sector consumes 77% of this value [56], *i.e.* 33.6 TWh. The Dutch Central Bureau for Statistics (CBS) reports service sector consumption of 30.6 TWh [57]. The service sector consumption in this paper amounts to 26.9 TWh for the entire Netherlands, *i.e.* 80% to 88% of the demand published by respectively PBL and CBS. The discrepancies in published data likely arise from the lack of unified definitions, an issue also raised by other researchers [10, 12, 25] and addressed further in Section 7. This quantitative validation indicates that the service sector profile estimation approach used in this paper can account for a substantial part of the Dutch service sector power demand. The remainder includes unaccounted for subsectors (*e.g.*, leisure), inaccuracies in subsector share estimations, and load profile deviations.

*4.3.2. Qualitative Validation*

The calculated service sector profiles are based on U.S. reference buildings. It remains an open question whether the use of U.S. buildings in the Dutch context causes deviations from real Dutch demand profiles. Perez-Lombard *et al.* [25] compared office energy end-use between U.S., Spain and the United Kingdom. End-use differences exist between the three countries. The differences between U.S. and the two European countries are however not larger than between the two European countries themselves. A similar qualitative conclusion can be drawn across the entire service sector by comparing the service sector end-use electricity consumption in the U.S. [58] and 29 European countries [16]. This suggests that using U.S. data for the Netherlands does not lead to larger errors than using data from another European country. Although undesirable, the practice of using data from other countries is currently common due to limited service sector data availability [10, 4].

Hereby it is important to note that the *shape* of the service sector demand profile, with a peak during the day, is similar across developed countries [10, 24, 25, 26]. It differs from the shape of household demand profiles, which typically peak in the evening [22, 59]. This observation qualitatively validates the use of U.S. profiles for the Dutch environment.

*4.4. Statistical Analysis*

Metric differences between residential loads-only and mixed loads are analysed for statistical significance using the two-sample t-test. Since multiple scenarios or categories are compared at once, the significance level is corrected using the Holm-Bonferroni correction to control the familywise error rate at 5%. For the first experiment, the correction is made for 121 comparisons. For second experiment, the correction is made for 150 comparisons.

## 5. Renewable Resource Integration Experiments

Two simulation experiments are carried out to study the impact of service sector loads on renewable resource integration. This impact is quantified using four metrics. The next paragraph outlines these metrics, the subsequent paragraphs provide details of the two experiments.

*5.1. Metrics*

The following renewable resource integration metrics are used in this paper:

- **Positive Mismatch.** Positive mismatch accounts for generation excess. It is calculated as the difference between generation and load when generation exceeds load.

- **Negative Mismatch.** Negative mismatch accounts for generation shortage. It is calculated as the difference between generation and load when load exceeds generation.

- **Renewable Energy Utilisation.** Renewable energy utilisation is the amount of renewable energy which can be used by the coinciding load. It is assumed that whenever renewable energy is available, it is utilised first. Only if no renewable energy is available, other (non-modelled) sources are used.

- **Self-Consumption.** Self-consumption is the ratio of renewable energy utilised by the coinciding loads and the total renewable energy generated.

*5.2. Experiment 1: Renewable Resource Penetration Scenarios*

Renewable resources considered in this paper are solar PV panels and wind turbines. For both solar PV and wind turbines, the installed generation capacity is varied between 0 MW and 525 MW with steps of 52.5 MW (121 scenarios in total). For the residential load case, 525 MW represents 300% of peak load (175 MW). For the mixed load case, 525 MW is 367% of peak load (143 MW), as mixed load has a flatter profile (see Fig. 1). The considered capacities are comparable to [8], where renewable resource capacity of up to 341% of peak load is considered for 2050.

In each scenario, the corresponding generation profile is calculated. This generation profile is combined with, on one hand, the demand profile of residential loads-only, and, on the other hand, with the demand profile of mixed loads. For each scenario and for each load type, a year-long hourly simulation is run. From the results, annual metrics are calculated and reported.

*5.3. Experiment 2: Time and Weather Dependency*

For a single scenario of solar PV and wind turbine penetration, this paper zooms in on the role time and weather conditions play in the impact the service sector has on renewable resource integration potential in urban systems with high renewable resource penetration. A novel time and weather dependency classification system is introduced to study the impact of different days of the week, times of the day, and weather conditions. The single scenario of solar PV and wind turbine penetration is obtained as a result of an area-constrained optimisation.

*5.3.1. Time and Weather Dependency Classification System*

In a power system with a high penetration of renewables, not only load variations, which mainly depend on the time of the day and the day of the week, determine the system state, but also weather variations, which govern renewable generation. To account for the future system dependency on both time and weather, this paper proposes a novel time and weather classification system. In this system, each hour of the year is classified according to four parameters: (1) day of the week, (2) time of day, (3) solar power generation and (4) wind power generation. Two categories are distinguished for the day of the week: weekday and weekend. Three categories are distinguished for the time of the day: night (00:00 - 08:00), day (08:00 - 16:00) and evening (16:00 - 00:00). Five categories are distinguished for both solar power generation and wind power generation. In both cases, the categories are based on quantiles. In total, 150 time and weather dependent categories are defined. Their frequency of occurrence is summarised in Table 2.

Table 2: Hours in 2014 classified according to the proposed time and weather dependency classification system. Numbers represent the hours per year in each category. The table distinguishes three parameters: time of the day, solar generation, and wind generation. The latter two are categorised using quantiles[1] and expressed as percentage of installed capacity. Note that for the sake of conciseness, this table makes no distinction between weekdays and weekends. The number of weekday and weekend hours in each category can be inferred statically: $5/7^{th}$ of all hours fall on weekdays and $2/7^{th}$ on weekends. For example, the number 627 in the first row, first column, represents the night hours in the reference year during which wind generation is between 0% and 5%, and solar generation is between 0% and 3%. The results corresponding to these hours are indicated in Figures 5 - 7 by a red arrow.

| Time of the Day | Wind Generation / Solar Generation | 0-5% | 5-13% | 13-26% | 26-71% | 71-100% | Total (hours) |
|---|---|---|---|---|---|---|---|
| **Night** | 0-3% | 627 | 419 | 426 | 409 | 414 | 2295 |
| | 3-9% | 92 | 53 | 40 | 29 | 29 | 243 |
| | 9-21% | 73 | 59 | 36 | 25 | 17 | 210 |
| | 21-40% | 35 | 48 | 37 | 13 | 9 | 142 |
| | 40-100% | 8 | 12 | 4 | 5 | 1 | 30 |
| **Day** | 0-3% | 31 | 26 | 48 | 68 | 99 | 272 |
| | 3-9% | 31 | 58 | 93 | 128 | 198 | 508 |
| | 9-21% | 50 | 78 | 106 | 157 | 173 | 564 |
| | 21-40% | 90 | 120 | 128 | 166 | 186 | 690 |
| | 40-100% | 132 | 226 | 177 | 196 | 155 | 886 |
| **Evening** | 0-3% | 540 | 526 | 532 | 450 | 419 | 2467 |
| | 3-9% | 21 | 58 | 45 | 39 | 18 | 181 |
| | 9-21% | 15 | 42 | 47 | 31 | 23 | 158 |
| | 21-40% | 7 | 24 | 28 | 30 | 10 | 99 |
| | 40-100% | 0 | 3 | 5 | 6 | 1 | 15 |
| | Total (hours) | 1752 | 1752 | 1752 | 1752 | 1752 | 8760 |

***Note 1.*** *During the night and in the evening more than 1752 hours fall in the first solar generation bracket (0-3%) because solar generation quantiles are calculated based on daylight hours.*

*5.3.2. Area-Constrained Renewable Mix Optimisation*

The single renewable resource penetration scenario is based on an area-constrained optimisation. The optimisation problem is formulated as a constrained multi-objective non-linear problem with design variables $x$ the number of solar PV panels and wind turbines: $x = [x_{PV}, x_{turbine}]$:

$$\begin{aligned}
\underset{x}{\text{minimize}} \quad & f(x) = p_{pos} \cdot M^+(x) + p_{neg} \cdot M^-(x) + p_{ren} \cdot R(x) \\
\text{subject to} \quad & 0 \leq x_{PV} \leq \phi \cdot A_{roof} \\
& 0 \leq x_{turbine} \leq (\phi - 1) \cdot A_{roof} \\
& x_{PV} + x_{turbine} \leq \phi \cdot A_{roof}
\end{aligned} \quad (7)$$

where:

$p_{pos}$ : weighting factor of positive mismatch ($p_{pos} > 0$, in this paper $p_{pos} = 1$)

$p_{neg}$ : weighting factor of negative mismatch ($p_{neg} > 0$, in this paper $p_{neg} = 1$)

$p_{ren}$ : weighting factor of renewable energy utilisation ($p_{ren} < 0$, in this paper $p_{ren} = -5$)

$A_{roof}$ : roof area available

$\phi$ : factor accounting for additional area available ($\phi \geq 1$, in this paper $\phi = 3$)

Positive mismatch $M^+(x)$, negative mismatch $M^-(x)$ and renewable energy utilisation $R(x)$ are all functions of the decision variables $x_{PV}$ and $x_{turbine}$ through their dependency on generation $G$. Generation $G$ is calculated as: $G = x_{PV} \cdot g_{PV} + x_{turbine} \cdot g_{turbine}$, with $g_{PV}$ and $g_{turbine}$ respectively the generation profiles of 1 m$^2$ PV and one 500 kW wind turbine.

This paper assumes that the total area available for renewable power generation is three times the size of the cumulative roof area of all the buildings considered (*i.e.* $\phi = 3$). This value represents 0.4% of Dutch land area, similar to the value originally estimated for the U.S. in [60]. The roof area itself is only available for solar power generation, while at most twice the roof area is available for wind power generation ($\phi - 1$ in Eq. 7). The footprint of a wind turbine is assumed to be 0.345 km$^2$/MW [61].

The problem at hand is a constrained multi-objective non-linear optimisation problem. The genetic algorithm in Matlab [27] is used to solve this problem. This algorithm relies on a population of possible individual solutions, which evolve to an optimal solution over a number of iterations. In each iteration, the best solutions are used to create solutions for the next iterations which are more likely to be close to the optimal solution. The algorithm terminates when the improvement in solutions falls below a threshold.

## 6. Case Study Results

This section presents the results of the two experiments conducted. The experiments quantify misestimations of renewable resource integration metrics that occur when the service sector is omitted, *i.e.* when mixed urban loads are represented by residential loads-only. The first experiment addresses misestimations in a broad range of renewable resource penetration scenarios. The second experiment zooms in on the misestimations on different days of the week, times of the day, and weather conditions for a single scenario.

### 6.1. Experiment 1: Renewable Resource Penetration Scenarios

Figure 4 shows annual average *differences* between (1) residential loads-only, and (2) mixed residential and service sector loads for four renewable resource integration metrics across a broad range of renewable resource penetration scenarios. Scenarios with solar and wind generation capacity of up to 525 MW are considered, *i.e.* 300% of peak load for the residential loads-only and 367% of peak load for the mixed loads (mixed loads have a flatter profile, see Figure 1).

#### 6.1.1. Mismatch

Figures 4a and 4b show respectively the annual average positive and negative mismatch *differences* between residential loads-only and mixed loads.

**Positive mismatch** represents renewable generation excess, *i.e.* renewable energy which cannot be used by the local loads. Positive mismatch *differences* indicate to what extent renewable generation excess is overestimated if residential loads-only are used instead of mixed loads. The positive mismatch difference is zero when solar and wind penetration equals zero, since no renewable power is generated. For all other penetration scenarios, differences increase with increasing solar penetration, while the variation as a function of wind is limited. Overall, Figure 4a shows that substituting mixed loads by residential loads-only leads to overestimation of generation excess. Results are statistically significant for solar penetration levels above 73% of peak load, and for wind penetration scenarios below 73% of peak load. Note that these cut-off values are based on the scenario step granularity of 36.5% of peak load.

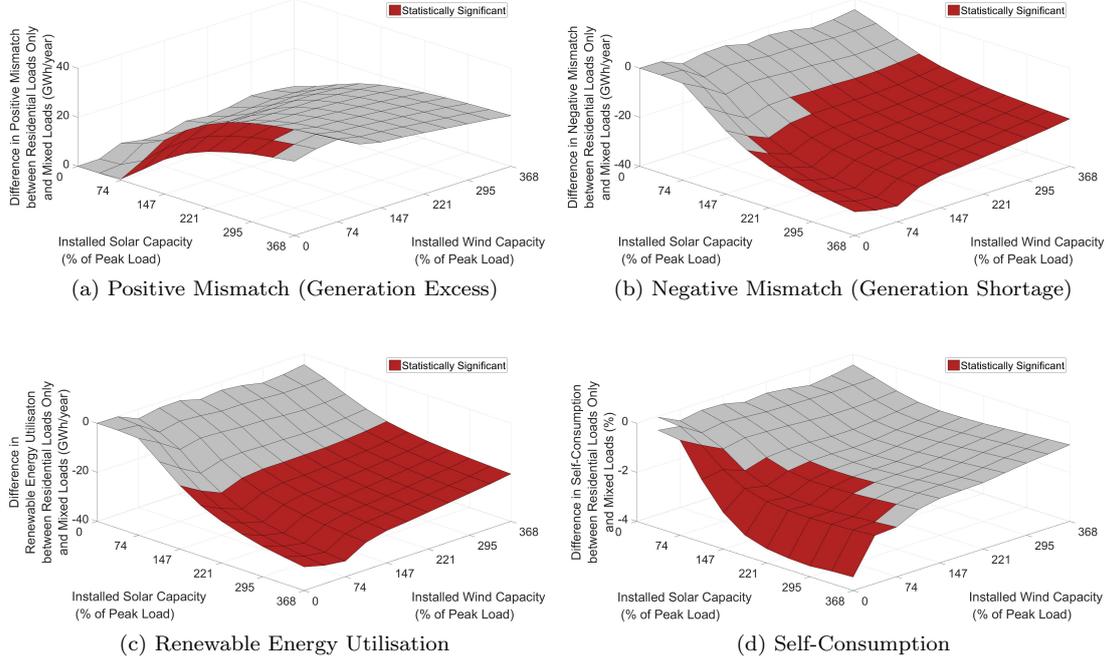

Figure 4: Annual average *differences* between residential loads-only and mixed residential and service sector loads. The XY plane represents scenarios of wind and solar penetration (expressed as percentage of peak load assuming mixed loads). Statistically significant differences are shown as red areas. Simulations for residential loads-only assume 203 005 households. Simulations for mixed loads assume 100 000 households and the corresponding number of service sector consumers as shown in Table 1. Note that the cumulative annual demand in both load types is equal (710 GWh/year).

**Negative mismatch** represents generation shortage, *i.e.* additional energy to be supplied by non-renewable resources. Negative mismatch *differences* indicate to what extent generation shortage is overestimated if residential loads-only are used instead of mixed loads. Negative mismatch difference is zero when solar and wind penetration equal zero as no renewable generation is available for either load type. Negative mismatch is larger in case of residential loads-only than in case of mixed loads, leading to negative mismatch differences below zero across all remaining scenarios. Overall, Figure 4b shows that substituting mixed loads by residential loads-only leads to overestimation of generation shortages. Results are statistically significant for scenarios with solar penetration above 110% of peak load (except for scenarios with very low installed wind capacity).

*6.1.2. Renewable Energy Utilisation*

Figure 4c shows renewable energy utilisation *differences* between residential loads-only and mixed loads. **Renewable energy utilisation** is the amount of renewable energy that is used by the coinciding demand. Renewable energy utilisation *differences* indicate to what extent the renewable energy utilisation is underestimated if residential loads-only are used instead of mixed loads. Renewable energy utilisation differences follow the same pattern as negative mismatch differences. For all scenarios, renewable energy utilisation is higher for the mixed loads than for the residential loads-only, thus the renewable energy utilisation difference is negative. Overall, Figure 4c shows that substituting mixed loads by residential loads-only leads to underestimation of renewable energy utilisation. Results are statistically significant for scenarios with solar capacity at or exceeding 147% of peak load, at any wind penetration.

*6.1.3. Self-Consumption*

Figure 4d shows self-consumption *differences* between residential loads-only and mixed loads. **Self-consumption** is the ratio of renewable energy utilised by the coinciding demand and the total renewable energy generated. Self-consumption *differences* indicate the extent to which the amount of generated renewable energy that can be used by the coinciding load is underestimated if residential loads-only are used instead of mixed loads. Self-consumption is highest when the penetration of renewable generation is low, it is undefined for zero penetration. If only a small amount of renewable power is generated, any type of coinciding load is sufficiently high to use it entirely. Self-consumption differences have a similar pattern as renewable energy utilisation differences, although differences at low wind penetration scenarios are more pronounced. Overall, Figure 4d shows that substituting mixed loads by residential loads-only leads to underestimation of self-consumption. Results are statistically significant for solar capacity scenarios above 73% of peak load and for wind penetration of at most 147% of peak load.

*6.1.4. Summary*

Differences in renewable resource integration metrics between residential loads-only and mixed loads are found across a broad range of scenarios. Considering all metrics together, statistically significant results are found in all scenarios except low solar. Note that non-significant results in low solar, low wind scenarios (left corner in Figures 4a-4d) occur because all metrics depend on the presence of renewable generation, which is very low in these scenarios. Overall, this experiment shows significant misestimations of annual average metrics if residential loads-only substitute mixed loads. The relative magnitude of these average annual differences is relatively small, up to approximately 5% of the total annual load. However, the differences between the metrics for residential loads-only and mixed loads vary throughout the year, depending on both time and weather conditions. These variations are assessed in the next experiment.

*6.2. Experiment 2: Time and Weather Dependency*

A power system with a high penetration of renewables is highly dependent on both time and weather. To study this dependency, all hours of the reference year (2014) are classified using the time and weather classification system introduced in this paper. Each category has four parameters: day of the week, time of the day, solar generation, and wind generation. In total, 150 time and weather dependent categories are analysed. An example of a category is: all weekday night (0:00 - 8:00) hours with solar generation between 0% and 3% of the installed capacity and wind generation between 0% and 5% of the installed capacity (this category is indicated on Figures 5, 6 and 7 with a red arrow). For each category, average metrics over all hours within that category are calculated and reported.

Results are shown in Figures 5, 6 and 7. In each figure, the upper row represents weekdays, the lower row – weekends. The columns represent three different times of the day: night, day and evening. Within each subfigure, 25 weather-dependent categories are shown. The three figures show the same metrics as considered for the renewable energy penetration scenarios, with positive and negative mismatches shown on one figure.

The results shown are obtained assuming an optimal renewable mix for the mixed loads: 399 MW solar PV and 30 MW wind turbines (*i.e.* total renewable capacity amounting to 245% of peak load in case of residential loads-only and 300% of peak load in case of mixed loads).

*6.2.1. Mismatch*

Figure 5 shows mismatch dependency on time and weather and compares residential loads-only and mixed loads. Positive mismatch indicates renewable generation excess. Negative mismatch indicates renewable generation shortage.

During weekdays and on weekend nights (Figure 5a-d), the mismatch is more positive for the residential loads-only than for the mixed loads. In the weekends, during the day and in the evening (Figure 5e-f), the mismatch is more positive for the mixed loads, although the differences are relatively small compared to the weekday categories. The largest differences occur on sunny weekdays (Figure 5a-c), and amount to up to 24% less mismatch between demand and supply in case of mixed loads than in case of residential loads-only.

The results obtained through the time and weather classification system can be used to identify critical combinations of time and weather. For instance, 62% of positive mismatches occur during weekdays at daytime when solar generation exceeds 40% of installed capacity, which corresponds to 7% of the time. Most negative mismatches (46%) occur during weekdays in the evening with solar generation below 3% of installed capacity, which corresponds to 20% of the time.

Statistical significance is not shown in the graph, yet is calculated as described in Section 4. Significant differences between mismatch results for the two load types are found for all data points on weekdays during the day (Figure 5b), as well as weekday and weekend evenings (Figure 5c and f) for low solar (generation below 3% of installed capacity). In other periods, statistically significant differences occur for some categories. The disparity in statistical significance between periods can be attributed to two factors: the number of data points and the relative difference between residential loads-only and mixed loads for a given period. First, as weather patterns are not dependent on the day of the week, weekdays have on average 2.5 times more data points per weather category than weekends. Second, during weekends and during night periods, the difference between residential loads-only and mixed loads is smaller than during other periods as most service sector activities are shut down.

*6.2.2. Renewable Energy Utilisation*

Figure 6 shows renewable energy utilisation dependency on time and weather and compares residential loads-only and mixed loads. Renewable energy utilisation is the amount of generated renewable energy that can be used by the coinciding loads. Higher renewable energy utilisation is better.

The differences in renewable energy utilisation between residential loads-only and mixed loads are most pronounced at high solar generation, both on weekdays and in weekends, and during all times of the day. Wind generation has limited effects as it represents only a small portion of the total renewable generation due to area constraints (see optimisation problem definition in Section 4). At higher solar generation levels, renewable energy utilisation is higher for the mixed loads than for the residential loads-only. The service sector consumption profile is more aligned with the solar power generation profile as both peak during the day. In the weekend, during day and evening hours (Figure 6e-f), the renewable energy utilisation at high solar irradiance levels is higher for residential loads-only. This corresponds to the fact that many service sector loads are minimal in the weekend. The largest differences between the two load cases occur on sunny weekdays (Figure 6a-c), and amount to up to 33% more renewable energy used directly by mixed loads than by residential loads-only.

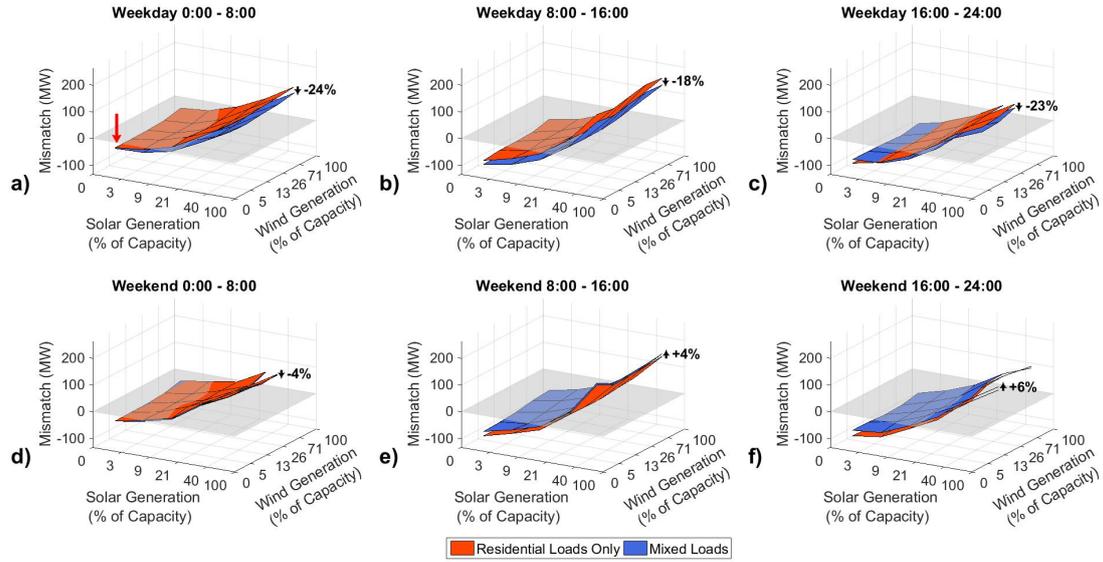

Figure 5: Mismatch dependency on time and weather: comparison between residential loads-only and mixed loads. The surface plots depict the average mismatches per time and weather category. The six subfigures represent different days of the week and times of the day. Within each subfigure, 25 weather-dependent categories are shown. Each category has four parameters: day of the week, time of the day, solar generation, and wind generation. For example, the category indicated by the red arrow represents all weekday night (0:00 - 8:00) hours of 2014 with solar generation between 0% and 3% of installed capacity and wind generation between 0% and 5% of installed capacity. The average mismatch in these hours is -41 MW for both load types. The black arrows on each subplot show the maximal relative changes in mismatch when mixed loads, instead of residential loads-only are taken into account (the values shown are for the most sunny hours). The values on the x- and y-axes are quantiles. Note that some high solar categories are missing because they do not occur in the modelled reference year. The results shown assume 399 MW solar PV and 30 MW wind turbines as installed capacity.

Statistically significant differences are found during weekdays at high solar generation levels for all periods (Figure 6a-c). Most renewable energy utilisation (26%) occurs during weekdays at daytime with high solar generation levels (above 40% of installed capacity), these categories correspond to 7% of the time. Further, 7% of the renewable energy is consumed during night and evening periods with lowest sun and highest wind (occurring 10% of the time).

### 6.2.3. Self-Consumption

Figure 7 shows self-consumption dependency on time and weather and compares residential loads-only and mixed loads. Self-consumption is the amount of renewable energy utilised relative to the amount generated.

As for the mismatch and renewable energy utilisation metrics, during weekdays (all periods) and on weekend nights the mixed loads performs better than the residential loads-only (Figure 7a-d). During weekend days and evenings the opposite is the case, although differences are again small (Figure 7e-f). As for other metrics, the largest differences are found on sunny weekdays (Figure 7a-c), mixed loads have a self-consumption of up to 32% higher than residential loads-only.

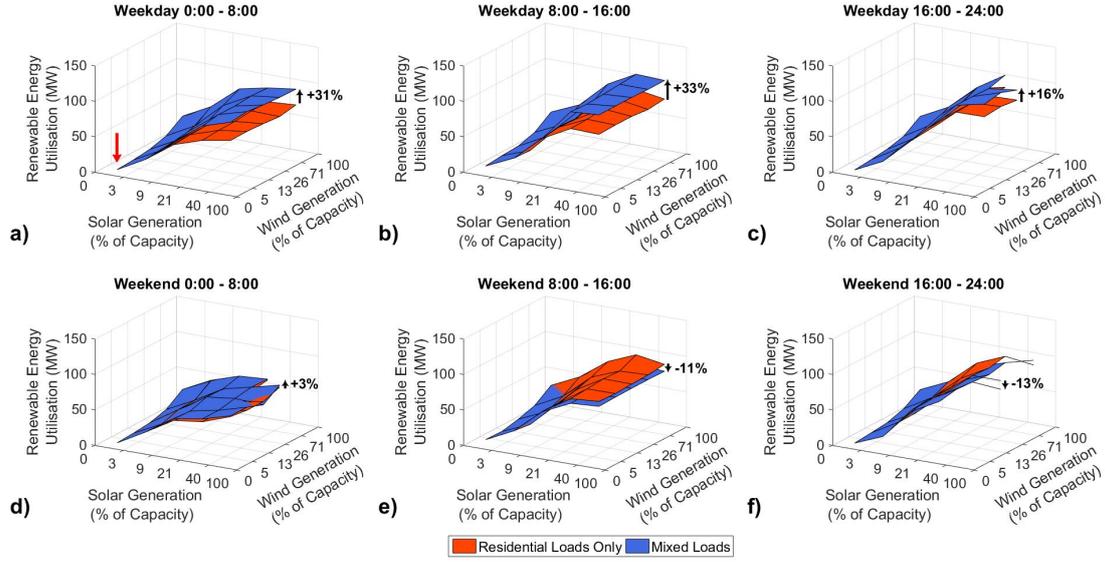

Figure 6: Renewable energy utilisation dependency on time and weather: comparison between residential loads-only and mixed loads. Surface plots depict the average renewable energy utilisation per time and weather category. The six subfigures represent different days of the week and times of the day. Within each subfigure, 25 weather-dependent categories are shown. The black arrows on each subplot show the maximal relative changes in renewable energy utilisation when mixed loads, instead of residential loads-only, are taken into account (the values shown are for the most sunny hours). The values on the x- and y-axes are quantiles. The red arrow indicates an example category introduced in Figure 5.

Statistically significant differences occur for similar categories as for mismatch. At low solar generation levels and at all wind generation levels, the self-consumption is 100%, meaning that all renewable power generated can be used by the modelled loads. As solar generation increases, self-consumption decreases. During weekdays the differences between the two load types are biggest (Figure 7b). In these periods the self-consumption decreases faster for the residential loads-only than for the mixed loads. This result illustrates that modelling only households underestimates the self-consumption of realistic mixed urban areas.

*6.2.4. Result Dependency on Load Assumptions in the Optimisation Step*

The results presented above rely on a renewable resource generation mix obtained by solving an optimisation problem assuming mixed loads. In this paper, the optimisation is constrained by area (see Section 5.3.2). This is the binding constraint for the number of wind turbines, regardless of the load type assumed. However, the optimal solar generation capacity changes with the load type. It is 15% lower if residential loads-only instead of mixed loads are assumed. The general trends for time and weather dependency as shown in Figures 5 - 7 remain similar if residential loads-only instead of mixed loads are assumed. However, overall mismatches become more negative, renewable energy utilisation decreases and self-consumption increases.

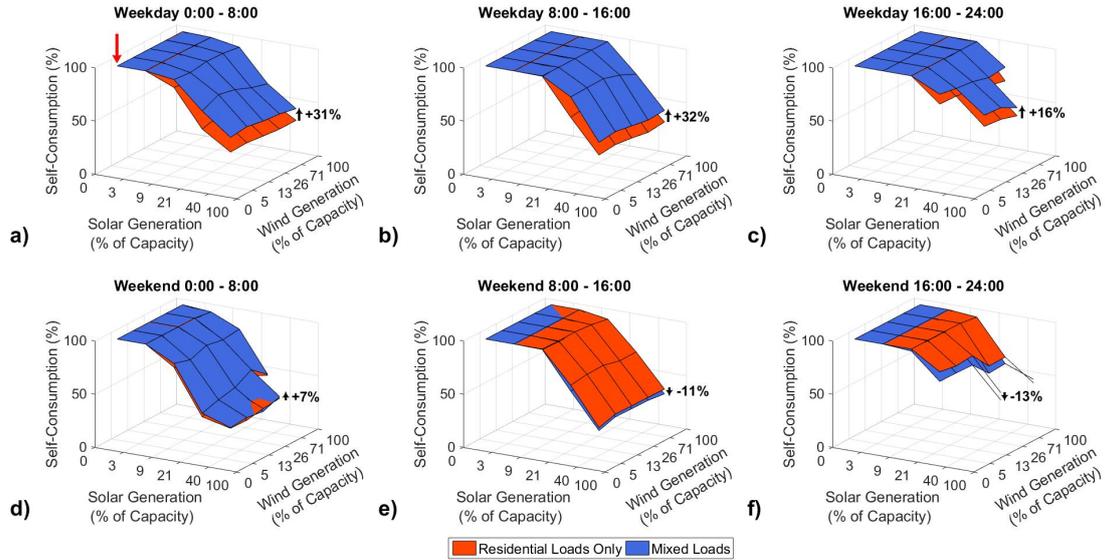

Figure 7: Self-consumption dependency on time and weather: comparison between residential loads-only and mixed loads. Surface plots depict the average self-consumption per time and weather category. The six subfigures represent different days of the week and times of the day. Within each subfigure, 25 weather-dependent categories are shown. The black arrows on each subplot show the maximal relative changes in self-consumption when mixed loads, instead of residential loads-only, are taken into account (the values shown are for the most sunny hours). The values on the x- and y-axes are quantiles. The red arrow indicates an example category introduced in Figure 5.

*6.2.5. Summary*

Renewable power integration metrics vary as a function of both time and weather. The results shown rely on the proposed time and weather classification system. Pronounced solar generation dependency is found for all metrics due to the high share of solar PV in the generation mix. Relative metric performance of residential loads-only and mixed loads differs per period. Overall, on weekdays (subplots a-c on Figures 5, 6 and 7) mixed loads lead to lower mismatches and higher renewable energy utilisation. During weekends (subplots d-f on Figures 5, 6 and 7) the contrary is the case. This difference can be attributed to service sector operation hours. Statistically significant differences between residential loads-only and mixed loads are primarily found on weekdays due to a larger number of datapoints per category and a larger difference between the two load type profiles. Overall, results show considerable differences (of up to 33%) between metrics calculated based on residential loads-only and those based on mixed loads.

## 7. Discussion

The integration of renewable energy resources in urban energy systems mandates a detailed understanding of the existing potential for renewable energy utilisation. In real urban areas, demand consists of a mix of residential and service sector loads. Existing urban energy system models primarily consider residential load profiles. As shown in this paper, omitting the service sector leads to misestimations of the potential for renewable energy utilisation in real urban areas. Currently, detailed measured demand data for the service sector are scarce. This paper (1) overcomes this lack of measured load profiles by devising synthetic service sector load profiles through a combination of a large number of different data sources, and (2) uses the obtained synthetic profiles to quantify misestimations of renewable resource integration metrics if service sector is not accounted for in urban areas. The four contributions of this paper are:

1. A systematic method to devise synthetic service sector load profiles

2. Quantification of renewable resource metric misestimations if the service sector is omitted, for a broad range of renewable resource penetration scenarios

3. A novel time and weather classification system

4. Quantification of renewable resource metric misestimations if the service sector is omitted, on different days of the week, times of the day, and weather conditions

Results reported in this paper can be valuable for researchers, practitioners, and decision-makers. More realistic urban demand profiles, based on both households and the service sector, can be used to extend urban energy system models, such as described by [1, 3, 4]. Decision-makers and practitioners can apply the reported results to improve grid planning, operation and management, for instance to guide interventions such as storage location, demand response programs, and grid reinforcement. The appropriate choice of such interventions depends on the timing and the extent of the mismatches between renewable generation and demand. Intervention choices based on misrepresented urban demand profiles (*e.g.*, profiles only accounting for residential loads in mixed urban areas), can lead to outcomes suboptimal for the real system.

This paper does not seek to determine which grid interventions are the most appropriate, as answering this question requires more detailed data than considered in this paper. Addressing this question is subject of further research. This paper focuses on quantitatively showing the overall importance of accounting for the service sector in the transition or urban areas to renewable generation. The next paragraphs discuss the four contributions of this paper.

*7.1. Synthetic Service Sector Load Profiles*

This paper models service sector demand using U.S. commercial reference building models and a combination of a large number of different U.S. and Dutch data sources. It proposes and implements a method to overcome the current lack of openly available, detailed measured service sector demand data for specific areas of interest. This method can considerably improve existing models of urban energy systems.

*7.1.1. Method Considerations*

The best validation of the proposed method relies on detailed, measured service sector profiles, the very issue this paper addresses. This is a chicken-or-egg problem. The proposed approach would not be necessary if detailed service sector or local (*e.g.*, municipality or neighbourhood-level) demand profiles are available. This is currently not the case.

The results of this paper show that the current approach to approximate urban demand by residential loads-only leads to statistically significant misestimations of renewable resource integration metrics. The method proposed in this paper can be used to estimate realistic urban demand profiles, and thus to improve estimations of renewable resource integration metrics.

In this paper the method is applied to an average Dutch urban area. The same approach can be applied to determine urban demand for other municipalities or neighbourhoods, provided sufficient local data are available. Neighbourhood-level urban demand modelling is currently researched by the authors.

*7.1.2. Result Generalisation*

This paper considers the Netherlands as a case study. It is an open question to what extent the results can *quantitatively* be generalised to other countries. The service sector composition and its share in the total national demand differ between countries [11, 12, 16]. This is not an issue in itself, the biggest challenge in comparing different regions arises due to inconsistencies in service sector definitions, as also underlined by other authors [10, 12, 25]. Even within a country, different sources provide different values for service sector power consumption (see Section 4). To improve service sector modelling, at least three issues need to be addressed: (1) inconsistent service sector definitions, (2) lack of openly available service sector data in general, and (3) lack of detailed (*e.g.*, hourly) service sector and local (*e.g.*, city or neighbourhood) load profiles in particular.

*Qualitatively*, our assumption is that the obtained results can be generalised to other developed countries because the *shape* of the service sector demand profile, with a peak during the day, is similar across developed countries [10, 24, 25, 26]. Based on the results presented in this paper, can be expected that the more important solar generation is in a country's renewable resource mix, the greater the impact of service sector loads is. Since solar power generation peaks during the day, it matches better with the service sector demand peak than with the household demand peak.

*7.2. Importance of the Service Sector across a Broad Range of Renewable Resource Penetration Scenarios (Experiment 1)*

Results obtained in experiment 1 show the impact of service sector loads on renewable resource integration across a broad range of renewable resource penetration scenarios. Statistically significant differences between renewable resource integration metrics for residential loads-only and mixed loads are found in all renewable resource penetration scenarios, except in those with high installed wind turbine capacity and low installed solar PV capacity, and those with few renewable resources (Figure 4).

Renewable generation scenarios with very high wind and low solar are highly unlikely due to physical constraints. Although the Netherlands currently produces ten times more renewable energy from wind (5300 GWh per year) than from solar (504 GWh per year) [62], this trend is unlikely to hold for high renewable resource penetrations. For equal installed capacity, wind turbines require considerably more area than solar panels. For instance, an installed renewable generation capacity of 367% of peak load would cover approximately 30% of the land area, or 20% of the off-shore Dutch Exclusive Economic Zone (up to 370 km off coast) [63] if wind turbines are used. If solar PV panels are used, the same installed capacity would cover only 1% of the land area.

From experiment 1 can be concluded that within the plausible range of scenarios, mixed loads lead to significantly less renewable energy excess, significantly less energy requirements from other non-renewable resources, and thus to a significantly higher renewable energy utilisation and significantly higher self-consumption. Although the future renewable mix is not known, these results show that service sector loads should be taken into account for renewable resource integration assessment in a broad range of plausible scenarios.

*7.3. Novel Time and Weather Classification System*

A renewable power system is highly dependent both on time and weather. Time governs diurnal, weekly and seasonal patterns in demand, and diurnal and seasonal patterns in solar generation patterns. Weather governs both solar and wind power generation, as well as some portion of the demand. Current power system metrics are assessed mainly from a time perspective [64]. The proposed time and weather classification system can contribute to the improvement of urban energy system models as it provides better insights on metric dependencies on time, weather and their interdependencies.

This paper proposes a novel time and weather dependency classification system which takes both time and weather into account. This classification system is flexible and can be readily applied to a wide range of dataseries. For the reference case used in this paper, categories are based on time intervals ne hour, full-year data, and five solar and wind energy generation categories (yielding 150 time and weather dependent categories). For other purposes, time interval, dataseries size and number of categories can be varied. For instance, the time and weather dependency classification system can be used with statistical data from multiple years to identify critical combinations of time and weather, to plan and manage distribution grid operations accordingly.

In the Results section such critical combinations are reported for the reference year 2014 (Section 6). The ability to identify such critical values as a function of time and weather and to assess their likelihood of occurrence is of importance for the design of grid interventions (*e.g.*, storage) and distribution grid management for power systems with a high share of renewable resources.

*7.4. Importance of the Service Sector in Specific Time and Weather Conditions (Experiment 2)*

Statistical analysis of the results obtained using the time and weather classification system shows significant metric differences between residential loads-only and mixed loads in a number of time and weather dependent categories. The most and largest differences (of up to 33%) are found on weekdays, in particular during sunny periods (over 1300 hours per year). These results demonstrate that using residential demand profiles to model mixed urban areas results in statistically significant metric misestimations. Such misestimations can have considerable impacts on, for instance, grid planning, operation and management choices.

The reported numerical results are based on the analysis of a single scenario. The following considerations indicate that the trends found can be generalised to other scenarios. First, significant annual differences in metrics are found across a broad range of scenarios (experiment 1). Second, the match of solar power generation with service sector power demand is better than with residential demand. Third, the mixed load profile is more constant than the residential profile, making it more likely that wind power generated at a random moment in time is used by mixed loads than by residential loads-only. Therefore, from the results obtained in experiment 2 can be generally concluded that during periods of high renewable power generation, the differences in metrics between residential loads-only and mixed loads are sufficiently large to necessitate the dedicated and detailed consideration of the service sector.

## 8. Conclusions and Future Work

This paper contributes to an improved understanding of future sustainable urban energy systems by showing the importance of including the service sector in energy system models. In the existing models, the service sector is often omitted due to the lack of detailed service sector load profiles for a specific area of interest, and the absence of a systematic method to devise them based on the very few available sources (such as [24]). This is the first systematic study addressing the impact of the service sector on renewable resource integration in urban areas. In this paper, a method is developed and implemented to devise synthetic service sector load profiles based on a combination of a large number of different openly available data sources. The obtained profiles are used to quantitatively show that omitting the service sector in urban energy systems leads to statistically significant misestimations of renewable resource integration metrics.

The proposed method and obtained results are being used for further research. Currently, the described method is extended by the authors to the neighbourhood level, to explore the local impact of storage. As residential and service sector loads are not evenly distributed in urban areas, concrete case studies of urban neighbourhoods are expected to provide further valuable local insights. Such insights are of importance for governments, distribution system operators, grid planners and new parties such as aggregators.

Future research directions include more extensive refinement and validation of the proposed method using measured service sector profile data, once they become available. Improving the proposed method further contributes to a better understanding of the measures needed to support the transition of cities to renewable resources.

## Acknowledgements

This work was supported by the Netherlands Organisation for Scientific Research (NWO) [grant number 408-13-012].

## References

[1] A. Fichera, M. Frasca, R. Volpe, Complex networks for the integration of distributed energy systems in urban areas, Applied Energy 193 (2017) 336–345. `doi:10.1016/j.apenergy.2017.02.065`.

[2] C. Hachem, Impact of neighborhood design on energy performance and GHG emissions, Applied Energy 177 (2016) 422–434. `doi:10.1016/j.apenergy.2016.05.117`.

[3] A. Alhamwi, W. Medjroubi, T. Vogt, C. Agert, GIS-based urban energy systems models and tools: Introducing a model for the optimisation of flexibilisation technologies in urban areas, Applied Energy 191 (2017) 1–9. `doi:10.1016/j.apenergy.2017.01.048`.

[4] J. Mikkola, P. D. Lund, Models for generating place and time dependent urban energy demand profiles, Applied Energy 130 (2014) 256–264. `doi:10.1016/j.apenergy.2014.05.039`.

[5] M. Paulus, F. Borggrefe, The potential of demand-side management in energy-intensive industries for electricity markets in Germany, Applied Energy 88 (2) (2011) 432–441. `doi:10.1016/j.apenergy.2010.03.017`.

[6] J. Torriti, M. G. Hassan, M. Leach, Demand response experience in Europe: Policies, programmes and implementation, Energy 35 (4) (2010) 1575–1583. `doi:10.1016/j.energy.2009.05.021`.

[7] P. Bertoldi, B. Atanasiu, Electricity Consumption and Efficiency Trends in the Enlarged European Union, Tech. rep., Europen Commission (2007).

[8] G. Brauner, W. D'Haseseleer, W. Gehrer, W. Glaunsinger, T. Krause, H. Kaul, M. Kleimaier, W. Kling, H. M. Prasser, I. Pyc, W. Schröppel, W. Skomudek, Electrical Power Vision 2040 for Europe, Tech. rep., EUREL, Brussels (2013).

[9] I. MacLeay, K. Harris, A. Annut, Digest of United Kingdom Energy Statistics 2014, Tech. rep., Department of Energy & Climate Change, London (2014).

[10] M. Jakob, S. Kallio, T. Bossmann, Generating electricity demand-side load profiles of the tertiary sector for selected European countries, in: 8th International Conference Improving Energy Efficiency in Commercial Buildings (IEECB14), no. April, 2014, p. 15. `doi:10.2790/32838`.

[11] Ofgem, Demand side response in the non-domestic sector, Tech. rep. (2012).

[12] N. Mairet, F. Decellas, Determinants of energy demand in the French service sector: A decomposition analysis, Energy Policy 37 (7) (2009) 2734–2744. `doi:10.1016/j.enpol.2009.03.002`.

[13] Energy Information Administration, Annual Energy Outlook 2015, Tech. rep., Energy Information Administration, Washington, D.C. (2013).

[14] D. Hostick, D. Belzer, S. Hadley, T. Markel, C. Marnay, M. Kintner-Meyer, Renewable Electricity Futures Study. End-use Electricity Demand, Tech. rep., National Renewable Energy Laboratory, Golden, CO (2012).

[15] P. Grünewald, J. Torriti, Demand response from the non-domestic sector: Early UK experiences and future opportunities, Energy Policy 61 (2013) 423–429. `doi:10.1016/j.enpol.2013.06.051`.

[16] T. Fleiter, S. Hirzel, M. Jakob, J. Barth, L. Quandt, F. Reitze, F. Toro, M. Wietschel, Electricity demand in the European service sector: A detailed bottom-up estimate by sector and by end-use, in: 6th International Conference on Improving Energy Efficiency in Commercial Buildings, IEECB Focus 2010, Frankfurt/Main, pp. 605–619.

[17] A. von Meier, Electric Power Systems: A Conceptual Introduction, John Wiley & Sons, 2006.

[18] E. Klaassen, J. Frunt, H. Slootweg, Assessing the Impact of Distributed Energy Resources on LV Grids Using Practical Measurements, in: 23rd International Conference on Electricity Distribution (CIRED), 2015, p. 5.

[19] I. Sartori, J. Ortiz, J. Salom, U. I. Dar, Estimation of load and generation peaks in residential neighbourhoods with BIPV: bottom-up simulations vs. Velander, in: World Sustainable Building 2014 Conference, 2014, pp. 17–24.

[20] USA Department of Energy, EnergyPlus Weather Data Sources, https://www.energyplus.net/weather/sources. Last accessed online 12-06-2017.

[21] Koninklijk Nederlands Meteorologisch Instituut, Uurgegevens van het weer in Nederland, https://www.knmi.nl/nederland-nu/klimatologie/uurgegevens. Last accessed online 12-06-2017.

[22] NEDU, Verbruiksprofielen, http://nedu.nl/portfolio/verbruiksprofielen. Last accessed online 12-06-2017.

[23] ENTSEO-E, Consumption Data, https://www.entsoe.eu/data/data-portal/consumption/Pages/default.aspx. Last accessed online 12-06-2017.

[24] M. Deru, K. Field, D. Studer, K. Benne, B. Griffith, P. Torcellini, B. Liu, M. Halverson, D. Winiarski, M. Rosenberg, M. Yazdanian, J. Huang, D. Crawley, U.S. Department of Energy commercial reference building models of the national building stock, Tech. rep., National Renewable Energy Laboratory (2011).

[25] L. Pérez-Lombard, J. Ortiz, C. Pout, A review on buildings energy consumption information, Energy and Buildings 40 (3) (2008) 394–398. `doi:10.1016/j.enbuild.2007.03.007`.

[26] G. Merei, J. Moshövel, D. Magnor, D. U. Sauer, Optimization of self-consumption and techno-economic analysis of PV-battery systems in commercial applications, Applied Energy 168 (2016) 171–178. `doi:10.1016/j.apenergy.2016.01.083`.

[27] MATLAB, Release 2015b, The MathWorks Inc., Natick, MA, USA, 2016.

[28] ECN, Energie-Nederland, Netbeheer Nederland, Energietrends 2014, Tech. rep. (2014).

[29] J. M. Bland, D. G. Altman, Regression towards the mean, British Medical Journal 308 (1994) 1499.

[30] American Hospital Association, Fast Facts on US Hospitals, http://www.aha.org/research/rc/stat-studies/fast-facts.shtml. Last accessed online 12-06-2017.

[31] Dutch Hospital Data, Kengetallen Nederlandse Ziekenhuizen 2014, Tech. rep., NVZ Vereniging van Ziekenhuizen in Nederland, Utrecht (2016).

[32] Centraal Bureau voor de Statistiek, Logiesaccommodaties; capaciteit, accommodaties, bedden, regio, http://statline.cbs.nl/StatWeb/publication/?VW=T&DM=SLNL&PA=82062NED&LA=NL. Last accessed online 12-06-2017.

[33] Invast Hotels, Hotelvastgoed Trendupdate, http://invasthotels.com/hotelvastgoed-trendupdate/2008/06/artikel-2/. Last accessed online 12-06-2017.

[34] R. Bak-Zeist, Kantoren in cijfers 2014. Statistiek van de Nederlandse kantorenmarkt, http://www.verbeek-bedrijfsmakelaars.nl/wp-content/uploads/2015/06/Kantoren-in-cijfers-2014.pdf. Last accessed online 12-06-2017.

[35] Compendium voor de Leefomgeving, Leegstand van kantoren, 1991-2015, http://www.clo.nl/indicatoren/nl2152-Leegstand-kantoren.html?i=36-177. Last accessed online 12-06-2017.

[36] W. Jiang, R. McBride, K. Jarnagin, M. Gowri, B. Liu, Technical Support Document: The Development of the Advanced Energy Design Guide for Highway Lodging Buildings, Tech. rep., Pacific Northwest National Laboratory (2008).

[37] Onderwijs in Cijfers, Aantal en omvang van instellingen in het primair onderwijs, http://www.onderwijsincijfers.nl/kengetallen/primair-onderwijs/instellingenpo/aantal-instellingen. Last accessed online 12-06-2017.

[38] Onderwijs in Cijfers, Aantal en omvang van vo-scholen, http://www.onderwijsincijfers.nl/kengetallen/voortgezet-onderwijs/instellingenvo/aantal-scholen. Last accessed online 12-06-2017.

[39] Compendium voor de Leefomgeving, Leegstand van winkels, 2004-2015, http://www.clo.nl/indicatoren/nl2151-Leegstand-winkels.html?i=36-177. Last accessed online 12-06-2017.

[40] Centraal Bureau voor de Statistiek, Elektriciteit in Nederland, Tech. rep., Den Haag (2015).

[41] DTZ Zadelhoff, De zekerheid van supermarkten, http://www.dtz.nl/media/90180/supermarkten-retail_folder_internet.pdf. Last accessed online 12-06-2017.

[42] Planbureau voor de Leefomgeving, Beleidsdossier detailhandel, http://themasites.pbl.nl/balansvandeleefomgeving/jaargang-2010/intensivering-verstedelijking-leefomgevingskwaliteit-en-woonwensen/bundeling-en-verdichting-verstedelijking/beleidsdossier-detailhandel. Last accessed online 12-06-2017.

[43] Rabobank, Rabobank Cijfers & Trends. Eetgelegenheden, https://www.rabobankcijfersentrends.nl/ index.cfm?action=branche&branche=Eetgelegenheden. Last accessed online 12-06-2017.

[44] Eurostat, Warehousing and transport support services statistics - NACE Rev. 2, http://ec.europa.eu/eurostat/statistics-explained/index.php/Archive:Warehousing_and_transport_support_services_statistics_-_NACE_Rev._2. Last accessed online 12-06-2017.

[45] USA Department of Energy, EnergyPlus Energy Simulation Software, https://energyplus.net/downloads. Last accessed online 12-06-2017.

[46] ASHRAE, ANSI/ASHRAE/IESNA Standard 90.1-2007. International Climate Zone Definitions, Tech. rep. (2008).

[47] G. Walker, Evaluating Mppt Converter Topologies Using a Matlab Pv Model, Journal of Electrical Electronics Engineering 21 (1) (2001) 49–56.

[48] Solarex, MSX-60 and MSX-64 Photovoltaic Modules, http://www.solarelectricsupply.com/media/custom/upload/Solarex-MSX64.pdf. Last accessed online 12-06-2017.

[49] Planbureau voor de Leefomgeving, Het potentieel van zonnestroom in de gebouwde omgeving van Nederland, Tech. rep. (2014).

[50] International Energy Agency, Technology Roadmap - Solar Photovoltaic Energy, Tech. rep. (2014).

[51] M. Z. Jacobson, M. A. Delucchi, Providing all global energy with wind, water, and solar power, Part I: Technologies, energy resources, quantities and areas of infrastructure, and materials, Energy Policy 39 (3) (2011) 1154–1169. doi:10.1016/j.enpol.2010.11.040.

[52] H. Lund, B. V. Mathiesen, Energy system analysis of 100 % renewable energy systems: The case of Denmark in years 2030 and 2050, Energy 34 (2009) 524–531. doi:10.1016/j.energy.2008.04.003.

[53] M. Patel, Wind and Solar Power Systems, CRC Press, 1999. doi:10.1201/9781420039924.

[54] EWT, High Yield 500kW Direct Drive Wind Turbine, http://www.ewtdirectwind.com/uploads/media/Leaflet.EWT500kW.pdf. Last accessed online 12-06-2017.

[55] Compendium voor de Leefomgeving, Aanbod en verbruik van elektriciteit, 1995-2015, http://www.clo.nl/indicatoren/nl0020-aanbod-en-verbruik-van-elektriciteit. Last accessed online 17-06-2017.

[56] ECN, Energie-Nederland, Netbeheer Nederland, Energietrends 2016, Tech. rep. (2016).

[57] Centraal Bureau voor de Statistiek, Aardgas en elektriciteit; leveringen openbaar net, bouw en dienstensector, http://statline.cbs.nl/Statweb/publication/?DM=SLNL&PA=82117ned&D1=1&D2=a&D3=a&VW=T. Last accessed 17-06-2017.

[58] L. Schwartz, M. Wei, W. Morrow, J. Deason, S. R. Schiller, G. Leventis, S. Smith, W. L. Leow, T. Levin, S. Plotkin, Y. Zhou, J. Teng, Electricity end uses , energy efficiency , and distributed energy resources baseline : Commercial Sector Chapter, Tech. Rep. January, Lawrence Berkeley National Laboratory (2017).

[59] J.-P. Zimmerman, M. Evans, J. Griggs, N. King, L. Harding, P. Roberts, C. Evans, Household Electricity Survey A study of domestic electrical product usage, Tech. rep., Intertek (2012).

[60] J. Turner, Renewable Energy: Generation, Storage, and Utilization, in: Carbon Management: Implications for R&D in the Chemical Sciences and Technology, 2001, pp. 111–126.

[61] P. Denholm, M. Hand, M. Jackson, S. Ong, Land-Use Requirements of Modern Wind Power Plants in the United States Land-Use Requirements of Modern Wind Power Plants in the United States, Tech. rep., National Renewable Energy Laboratory (2009).

[62] Centraal Bureau voor de Statistiek, Hernieuwbare energie in Nederland, Tech. rep., Den Haag (2013).

[63] Marineregions, Netherlands - MGGID 5568, http://www.marineregions.org/eezdetails.php?mrgid=5668. Last accessed online 12-06-2017.

[64] H. Hahn, S. Meyer-Nieberg, S. Pickl, Electric load forecasting methods: Tools for decision making, European Journal of Operational Research 199 (3) (2009) 902–907. `doi:10.1016/j.ejor.2009.01.062`.

[65] Compendium voor de Leefomgeving, Bevolkingsomvang en aantal huishoudens 1980-2015, http://www.clo.nl/indicatoren/nl0001-Bevolkingsomvang-en-huishoudens.html?i=15-1. Last accessed online 12-06-2017.

# Appendix

## A1. Calculation of Service Sector Scaling Factors

United States Department of Energy (U.S. DOE) commercial reference building data [24] are used to model service sector electricity demand in the Netherlands. Each U.S. reference building is *scaled* to represent its Dutch counterpart. Scaling factors are calculated based on Dutch national data, or data obtained from subsector-specific organisations (see further). The calculations are described below per U.S. reference building (group). Where-ever available, data for the reference year 2014 are used to ensure the best correspondence with other data used in this study. Otherwise, the latest available data are used.

This study uses a reference urban environment of 100 000 households. To obtain the number of reference buildings representative for this urban environment, the number of reference buildings obtained for the Netherlands is divided by 75.9, as the total number of households in the Netherlands is 7.59 million [65]. The final result is rounded to the nearest integer.

### A1.1. Hospitals

The scaling factors for U.S. reference building "Hospital" is calculated based on the number of patient beds. The average number of beds in hospitals in the U.S. is 161 [30]. The average number of beds in hospitals in the Netherlands is 316 [31]. The scaling factor can thus be calculated as $316/161 = 1.96$. There are 134 hospitals in the Netherlands [31], which can be represented by $134 * 1.96 = 263$ U.S. reference buildings of the type "Hospital". In an urban environment of 100 000 households, Dutch hospitals can be represented by $263/75.9 = 3$ U.S. reference buildings of the type "Hospital".

### A1.2. Hotels

The scaling factors for U.S. reference buildings "Large Hotel" and "Small Hotel" are calculated based on the number of hotel rooms. The number of rooms in the "Large Hotel" U.S. reference building is 300 [24]. The number of rooms in the "Small Hotel" U.S. reference building is 77 [24]. There are 3185 hotels in the Netherlands with a total of 112 565 rooms [32]. 2.6% of the hotels (83 hotels) in the Netherlands have more than 200 rooms per hotel [33]. These hotels are represented by the "Large Hotel" reference building. Hotels with less than 200 rooms are represented by the "Small Hotel" reference building. Large Dutch hotels are assumed to have an average of 250 rooms per hotel, and can thus be represented by $83 * 250/300 = 69$ U.S. reference buildings of the type "Large Hotel". The remaining 91 865 rooms can be represented by 1193 U.S. reference buildings of the type "Small Hotel". In an urban environment of 100 000 households, Dutch hotels can be represented by $83/75.9 = 1$ U.S. reference buildings of the type "Large Hotel" and $1193/75.9 = 16$ U.S. reference buildings of the type "Small Hotel".

### A1.3. Offices

The scaling factors for U.S. reference buildings "Large Office", "Medium Office" and "Small Office" are calculated based on office floor area. The floor area of the "Large Office" U.S. reference building is 46 320 m$^2$ [24]. The floor area of the "Medium Office" U.S. reference building is 4982 m$^2$ [24]. The floor area of the "Small Office" U.S. reference building is 511 m$^2$ [24]. Office floor area distribution for the Netherlands is shown in Table A1 [34]. The total (used) office floor area in the Netherlands is 49.55 million m$^2$ [35].

Dutch offices with an area larger than 10 000 m$^2$ are represented by $31\ 263\ 884/46\ 320 = 675$ "Large Office" U.S. reference buildings. Dutch offices with an area between 1000 and

Table A1: Floor area distribution of offices in the Netherlands [34]. The average area is obtained from the minimum and maximum areas. The number of offices is the solution of a set of equations which satisfy both the third column and the total (used) office floor area in the Netherlands, which equals 49.55 million m$^2$ [35]. The last column illustrates this.

| Input Data [34] | | | Calculated Values | | |
|---|---|---|---|---|---|
| Min. Area | Max. Area | Share | Average Area | Number of Offices | Total Area |
| (m$^2$) | (m$^2$) | (%) | (m$^2$) | (-) | (m$^2$) |
| 500 | 1000 | 5 | 750 | 326 | 244 249 |
| 1000 | 2500 | 19 | 1750 | 1238 | 2 165 675 |
| 2500 | 5000 | 23 | 3750 | 1498 | 5 617 729 |
| 5000 | 10000 | 21 | 7500 | 1368 | 10 258 462 |
| 10000 | (-) | 32 | 15000 | 2048 | 31 263 884 |

10 000 m$^2$ are represented by (2 165 675 + 5 617 729 + 10 258 462)/4982 = 3569 "Medium Office" U.S. reference buildings. Dutch offices with an area between 500 and 1000 m$^2$ are represented by 244 249/511 = 478 "Small Office" U.S. reference buildings.

In an urban environment of 100 000 households, Dutch offices can be represented by 675/75.9 = 9 U.S. reference buildings of the type "Large Office", 3569/75.9 = 47 U.S. reference buildings of the type "Medium Office" and 478/75.9 = 6 U.S. reference buildings of the type "Small Office".

*A1.4. Schools*

The scaling factors for U.S. reference buildings "Primary School" and "Secondary School" are calculated based on student numbers. The number of students in the "Primary School" U.S. reference building is 650 [24, 36]. The number of students in the "Secondary School" U.S. reference building is 1200 [24, 36]. The number of primary schools in the Netherlands is 7155, with an average of 219 students per school [37]. These schools can be represented by 7155 ∗ 219/650 = 2411 U.S. reference buildings of the type "Primary School". The number of secondary schools in the Netherlands is 642 with an average of 1295 students per school [38]. These schools can be represented by 642 ∗ 1295/1200 = 693 U.S. reference buildings of the type "Secondary School". In an urban environment of 100 000 households, Dutch schools can be represented by 2411/75.9 = 32 U.S. reference buildings of the type "Primary School" and 693/75.9 = 9 U.S. reference buildings of the type "Secondary School".

*A1.5. Retail*

The scaling factor for U.S. reference building "Stand Alone Retail" is calculated based on floor area. The floor area of "Stand Alone Retail" U.S. reference building is 2294 m$^2$ [24]. The total (used) retail area in the Netherlands is 30 775 168 m$^2$ [39]. Thus, retail in the Netherlands can be represented by 30 775 168/2294 = 13 416 U.S. reference buildings of the type "Stand Alone Retail". In an urban environment of 100 000 households, Dutch retail can be represented by 13 416/75.9 = 177 U.S. reference buildings of the type "Stand Alone Retail".

*A1.6. Supermarkets*

The scaling factor for U.S. reference building "Supermarket" is calculated based on floor area. The floor area of "Supermarket" U.S. reference building is 4181 m$^2$[24]. The total supermarket area in the Netherlands is 3 781 699 m$^2$[41] (cross-referenced with [40, 42]). Thus, supermarkets in the Netherlands can be represented by 3 781 699/4181 = 904 U.S. reference buildings of the type "Supermarket". In an urban environment of 100 000 households, Dutch retail can be represented by 904/75.9 = 12 U.S. reference buildings of the type "Supermarket".

*A1.7. Restaurants*

The scaling factors for U.S. reference buildings "Full Service Restaurant" and "Quick Service Restaurant" are calculated based on the number of restaurants in the Netherlands. The number of restaurants of different types are obtained from [43] and are summarised in Table A2. The table indicates which restaurants are consideblack to be equivalent with "Full Service Restaurant" and which with "Quick Service Restaurant". The distinction is made based on the expected time customers spend at a restaurant. The total number of restaurants in the Netherlands which can be represented by U.S. reference building of the type "Full Service Restaurant" is 12 903. The total number of restaurants in the Netherlands which can be represented by U.S. reference building of the type "Quick Service Restaurant" is 14 372. In an urban environment of 100 000 households, Dutch restaurants can be represented by 12 903/75.9 = 170 U.S. reference buildings of the type "Full Service Restaurant" and 14 372/75.9 = 190 U.S. reference buildings of the type "Quick Service Restaurant".

*A1.8. Warehouses*

The scaling factors for U.S. reference building "Warehouses" are calculated based on the power consumption and the number of employees in warehousing in the Netherlands and the United Kingdom (UK) due to lack of Dutch data. The annual electricity consumption of "Warehouse" U.S. reference building is 239 MWh per year [24]. The total warehouses electricity consumption in the UK is 12 TWh per year [11]. The total number of employees in warehousing in the UK is 334.1 thousand. The total number of employees in warehousing in the Netherlands is 82.6 thousand [44]. Warehouses in the Netherlands are therefore estimated to be represented by 12 000/0.239∗826 000/3 341 000 = 12 397 U.S. reference buildings of the type "Warehouse". In an urban environment of 100 000 households, Dutch warehouses can be represented by 12 397/75.9 = 164 U.S. reference buildings of the type "Warehouse".

Table A2: Number of restaurants of different types in the Netherlands [43].

| Restaurant Type (Dutch Name) | Restaurant Type (English Translation) | Number | U.S. Restaurant Equivalent |
| --- | --- | --- | --- |
| Lunchroom/Tearoom | Lunchroom/Tearoom | 1929 | Full |
| Selfservice restaurant | Self-service Restaurant | 162 | Quick |
| Broodjeszaak/Croissanterie | Sandwichbar/Eat-in Bakery | 937 | Quick |
| Internationale fastfoodketen | International Fast Food Chain | 325 | Quick |
| Cafetaria/Snackbar | Cafeteria/Snack Bar | 5888 | Quick |
| Shoarma/Grillroom/Kebab | Shoarma/Grillroom/Kebab | 1926 | Quick |
| IJssalon | Ice Cream Parlor | 876 | Quick |
| Restaurant (buitenlandse keuken) | Restaurant (International Cuisine) | 5667 | Full |
| Restaurant (nationale keuken) | Restaurant (National Cuisine) | 4663 | Full |
| Pizzeria | Pizzeria | 644 | Full |
| Eetcafe | Pub | 3318 | Quick |
| Creperie/Pannenkoeken | Creperie/Pancakes | 429 | Quick |
| Recreatiemeer-/Strandpaviljoen | Recreation Park/Beach Pavilion | 511 | Quick |